\documentclass[letterpaper, 10 pt, conference]{ieeeconf}  

\IEEEoverridecommandlockouts                              

\usepackage{algorithm2e}
\usepackage{algorithmicx}

\usepackage{tikz}
 \usetikzlibrary{arrows,shapes,positioning}
 \usetikzlibrary{decorations.markings}
 \tikzstyle arrowstyle=[scale=1]
 \tikzstyle directed=[postaction={decorate,decoration={markings,mark=at position .65 with {\arrow[arrowstyle]{stealth}}}}]
 \tikzstyle reverse directed=[postaction={decorate,decoration={markings,mark=at position .65 with {\arrowreversed[arrowstyle]{stealth};}}}]
\usepackage{amsmath}
\usepackage{epstopdf}
\usepackage{latexsym, amssymb,booktabs}
\usepackage{graphicx}
\usepackage{setspace}
\usepackage{amsmath, amsbsy}
\usepackage{amsopn, amstext}
\usepackage{hyperref,url}
\usepackage{floatrow}
\newfloatcommand{capbtabbox}{table}[][\FBwidth]
\usetikzlibrary{arrows}

\DeclareMathOperator{\Col}{Col}

\DeclareGraphicsExtensions{.eps,.pdf}

\def\ra{\rightarrow}

\def\d{\delta}

\def\D{\Delta}

\newcommand{\R}{{\mathbb R}}

\newtheorem{thm}{Theorem}[section]

\newtheorem{dfn}[thm]{Definition}

\newtheorem{exa}[thm]{Example}

\title{A Note On Orthogonal Decomposition of Finite Games}
\author{Changxi Li$^{1}$,  Daizhan Cheng$^{2}$, Fenghua He$^{1}$,
\thanks{*This work is supported partly by the National Natural Science Foundation of China (NSFC) under Grants 61473099, 61773371, 61733018 and 61333001.}
\thanks{$^{1}$ Changxi Li and Fenghua He are with Harbin Institute of Technology, Harbin 150001, P.~R.~China {\tt\small  changxileehit@gmail.com, hefenghua@hit.edu.cn}}
\thanks{$^{2}$ Daizhan Cheng is with Key Laboratory of Systems and Control,
		Academy of Mathematics and Systems Sciences, Chinese Academy of Sciences,
		Beijing 100190, P. R. China {\tt\small   dcheng@iss.ac.cn}}
\thanks{ Corresponding author: Changxi Li. }
}

\begin{document}

	\maketitle
	\thispagestyle{empty}
	\pagestyle{empty}


\begin{abstract}
Various decomposition of finite games have been proposed. The inner product of vectors plays a key role in the decomposition of finite games. This paper considers the effect of different inner products on the  orthogonal decomposition of finite games.  We find that only when the compatible condition is satisfied,  a common decomposition  can be induced by the standard inner product and the weighted inner product. To explain the result, we studied the existing decompositions, including potential based decomposition, zero-sum based decomposition, and symmetry based decomposition.
\end{abstract}


\section{Introduction}

Due to the widespread applications of game theory, many researchers begin studying the topological structure of finite games. Decomposition of finite games is the main technique \cite{whs10}-\cite{gsz17}. 

Different decompositions are proposed from different point of views. i) Using Helmholtz decomposition theorem, \cite{can11} proposed a \emph{potential based decomposition} of finite games, where the space of finite games ${\cal G}_{[n;k_1,\cdots,k_n]}$ is decomposed into a canonical  sum of the pure potential subspace ${\cal P}$, the non-strategic subspace ${\cal N}$, and the pure harmonic subspace ${\cal H}$,
\begin{align*}
{\cal G}_{[n;k_1,\cdots,k_n]}=\rlap{$\underbrace{\phantom{\quad{\cal P}\quad\oplus\quad{\cal N}}}_{Potential\quad games}$}\quad{\cal P}\quad\oplus\quad
\overbrace{{\cal N}\quad\oplus\quad{\cal H}}^{Harmonic\quad games}.
\end{align*}
ii) \cite{che16} reinvestigated the \emph{potential based decomposition} on  Euclidean space, and bases of all subspaces in the decomposition are provided. iii) \cite{shh16} proposed three kinds of decompositions of finite games: \emph{zero-sum based decomposition}
\begin{align*}
{\cal G}_{[n;k_1,\cdots,k_n]}={\cal Z}~\oplus~{\cal C},
\end{align*}
\emph{normalization based  decomposition}
\begin{align*}
{\cal G}_{[n;k_1,\cdots,k_n]}={\cal L}~\oplus~{\cal E},
\end{align*}
and \emph{zero-sum equivalent potential based decomposition}
\begin{align*}
{\cal G}_{[n;k_1,\cdots,k_n]}=({\cal L}\cap{\cal C})~\oplus~{\cal B}~\oplus~({\cal L}\cap{\cal Z}),
\end{align*}
where ${\cal Z}, {\cal C}, {\cal L}, {\cal E}$ and ${\cal B}=({\cal Z}+{\cal E})\cap({\cal C}+{\cal E})$ are the subspace of zero-sum games, common interest games, normalized games, non-strategic games, and zero-sum equivalent potential games, respectively. iv) \emph{Symmetry based decomposition} of finite games are proposed in \cite{cx18,dct18}, which is shown as follows
\begin{align*}
{\cal G}_{[n;l,\cdots,l]}=\mathcal{S}\oplus\mathcal{K},
\end{align*}
where $\mathcal{S}$ is the subspace of the symmetry games and $\mathcal{K}$ is the orthogonal complement of $\mathcal{S}$.

As far as we know,  different technical tools are used for various decompositions. For example, the technical tools used in \cite{can11} and \cite{che16} are Helmholtz decomposition theorem and semi-tensor product (STP) of matrices respectively, where different inner products for the space of games are defined. But surprisingly, \cite{can11} and \cite{che16} have provided the same decomposition (\emph{potential based decomposition}) independently. Why can different inner products  induce the same decomposition? Contingency or necessity? If it is a contingency, then under what conditions will the decomposition be the same using Helmholtz decomposition theorem and STP? This note aims at answering all the questions.

The rest of this paper is organized as follows:  Section 2 provides some preliminaries including finite games theory, matrix expression of finite games. Section 3 investigates the effect of different inner products on the orthogonal decomposition of finite games.  Section 4 is a brief conclusion.
For statement ease, we give some notations:
\begin{enumerate}
\item ${\cal M}_{m\times n}$: the set of $m\times n$ real matrices.
\item ${\cal D}:=\left\{0,1\right\}$.
\item $\d_n^i$: the $i$-th column of the identity matrix $I_n$.
\item $\D_n:=\left\{\d_n^i\vert i=1,\cdots,n\right\}$.
\item ${\bf 1}_{\ell}=(\underbrace{1,1,\cdots,1}_{\ell})^T$.
\item ${\bf 0}_{p\times q}$: a $p\times q$ matrix with zero entries.
\item A matrix $L\in {\cal M}_{m\times n}$ is called a logical matrix if the columns of $L$ are of the form of $\d_m^k$. That is, $\Col(L)\subset \D_m$.
      Denote by ${\cal L}_{m\times n}$ the set of $m\times n$ logical matrixes.
\item If $L\in {\cal L}_{n\times r}$, by definition it can be expressed as $L=[\d_n^{i_1},\d_n^{i_2},\cdots,\d_n^{i_r}]$. For the sake of compactness, it is briefly denoted as
      $L=\d_n[i_1,i_2,\cdots,i_r]$.
\item Let $N=\{1,2,\cdots,n\}$. ${\bf S}_n$ is the set of permutations of elements of $N$.
\end{enumerate}
\section{Preliminaries}
This section gives a very brief review on semi-tensor product of matrices and game theory. Plese refer to \cite{gib92}, \cite{che12} for more details.

\subsection{Finite Games}
\begin{dfn}\label{d1.1} A finite non-cooperative game $G=\{N,S,c\}$ is a triple, where
\begin{enumerate}
\item[(a)] $N=\{1,2,\cdots,n\}$ is the set of players;
\item[(b)] $S=\prod_{i=1}^nS_i$ is the strategy profile of the game, with $S_i=\{1,2,\cdots,k_i\}$ as the set of strategies of player $i\in N$;
\item[(c)] $c=\{c_1,\ldots,c_n\}$ is the set of payoff functions, where $c_i:S\ra \R$ is the payoff function of player $i\in N$.
\end{enumerate}
\end{dfn}

Denote by ${\cal{G}}_{[n;k_1,\cdots,k_n]}$ the set of finite games with $|N|=n$  and $|S_i|=k_i, \forall i=\in N$.
There are many classes of games. The following are the definitions of some special games used in this note.
\begin{itemize}
\item \emph{Zero-sum games}: a finite non-cooperative game $G\in{\cal{G}}_{[n;k_1,\cdots,k_n]}$ is a zero-sum game, if and only if,
$$\sum_{i=1}^n c_i(s)=0,\ \forall s\in S.$$
\item \emph{Common interest games}: a finite non-cooperative game $G\in{\cal{G}}_{[n;k_1,\cdots,k_n]}$ is a common interest game, if and only if,
$$ c_i(s)=c_j(s), \ \forall s\in S, \forall i,j\in N.$$
\item \emph{Normalized  games}: a finite non-cooperative game $G\in{\cal{G}}_{[n;k_1,\cdots,k_n]}$  is a normalized game, if and only if,
$$ \sum_{x_i\in S_i} c_i(x_i, s_{-i})=0,\ \forall s_{-i}\in S_{-i},\forall i\in N,$$
where $S_{-i}=\prod_{j=1,j\neq i}^nS_j$.
\item \emph{Potential  games}: a finite non-cooperative game $G\in{\cal{G}}_{[n;k_1,\cdots,k_n]}$ is potential if and only if there is a function $P(s)$, called the potential function, such that
\begin{align*}
c_i(x_i,s_{-i})-c_i(y_i,s_{-i})=P(x_i,s_{-i})-P(y_i,s_{-i}) \\
\forall x_i,y_i\in S_i,\;\forall s_{-i}\in S_{-i},\; i=1, \cdots, n.
\end{align*}
\item \emph{Harmonic games}: a finite non-cooperative game $G\in{\cal{G}}_{[n;k_1,\cdots,k_n]}$ is a harmonic game, if and only if, it is a zero-sum game and a normalized game, i.e.,
$$\sum_{i=1}^n c_i(s)=0,\ \forall s\in S,$$
$$ \sum_{x_i\in S_i} c_i(x_i, s_{-i})=0,\ \forall s_{-i}\in S_{-i},\forall i\in N.$$
\item \emph{Non-strategic games}: a finite non-cooperative game $G\in{\cal{G}}_{[n;k_1,\cdots,k_n]}$ is a non-strategic game, if and only if, for any $i\in N$ that
\begin{align*}
c_i(x_i, s^{-i})-c_i(y_i,s^{-i})=0,\quad \forall x_i,y_i\in S_i,\;\forall s^{-i}\in S^{-i}.
\end{align*}
\item \emph{Symmetry games}: a finite non-cooperative game $G\in{\cal{G}}_{[n;l,\cdots,l]}$ is a symmetry game, if and only if,   for any $\sigma\in {\bf S}_n$
\begin{align*}
c_i(s_1,\cdots,s_n)=c_{\sigma(i)}(s_{\sigma^{-1}(1)}, s_{\sigma^{-1}(2)},\cdots,s_{\sigma^{-1}(n)}),\\
 s_i\in S_i, i=1,\cdots,n.
\end{align*}
\end{itemize}

\subsection{Matrix Expression of Finite Games}
The tool used in this paper is  the semi-tensor product (STP) of matrices \cite{che12}, which is a generalization of conventional matrix product.
\begin{dfn}\label{d1.6} Let $A\in {\cal M}_{m\times n}$ and  $B\in {\cal M}_{p\times q}$ and $t=lcm(n,~p)$ be the least common multiple of $n$ and $p$. The semi-tensor product (STP) of $A$ and $B$ is defined as
\begin{align*}
A\ltimes B:=\left(A\otimes I_{t/n}\right)\left(B\otimes I_{t/p}\right)\in {\cal M}_{mt/n\times qt/p}.
\end{align*}
\end{dfn}

Identify $j\sim \d_{k_i}^j$, which is called the vector expression to strategies $j\in S_i$.
Using STP and the vector expression to strategies $s_i\in S_i$, $i=1,\cdots,n$, the strategy profile $s=\prod_{i=1}^n s_i$ can be expressed as
$$
s=\ltimes_{i=1}^ns_i.
$$
Under this expression, each payoff function $c_i$ becomes a mapping $c_i:\D_k\ra \R$, where $k=\prod_{i=1}^nk_i$. Hence for each $c_i$ we can find a unique row vector $V_i\in \R^k$ such that
\begin{align*}
\begin{array}{ccl}
c_i(s_1,\cdots,s_n)=V_ix, \quad i=1,\cdots,n.
\end{array}
\end{align*}
$V_i$ is called the structure vector of $c_i(x)$. A finite game $G\in {\cal G}_{[n;k_1,\cdots,k_n]}$ is uniquely determined by $\{V_i|i=1,\cdots,n\}$.
Denote the payoff vector by $V_G:=[V_1,\cdots,V_n]\in \R^{nk}$.
Then it is clear that ${\cal G}_{[n;k_1,\cdots,k_n]}$ has a natural vector space structure as $\R^{nk}$.

To  illustrate the  vector expression of finite games, we provide the following example.
\begin{exa}\label{e2.1}
Consider a three-player  game $G$. Each player has two strategies $S_i=\{1,2\},i=1,2,3$, and the payoffs of $G$ are  described as in Table \ref{Tab3.3}.
\begin{table}[!htbp]
  \centering
  \caption{Payoff Matrix of Example \ref{e2.1}}\label{Tab3.3}
  \begin{tabular}{cccccccccc}
    \hline $V_i$$\backslash\mbox{s}$ &$111$&$112$&$121$&$122$&$211$&$212$&$221$&$222$\\
    \hline $V_1$ & $26$ & $9$ & $12$ & $4$  & $14$ & $6$ & $14$ & $6$\\
           $V_2$ & $-5$ & $-5$ & $2$ & $2$  & $2$ & $2$ & $4$ & $4$\\
           $V_3$ & $18$ & $10$ & $4$ & $5$  & $7$ & $8$ & $7$ & $8$\\
    \hline
  \end{tabular}
\end{table}

The payoff vector of $G$ is
\begin{align*}
\begin{array}{ccl}
V_G&=&[V_1,V_2,V_3]\\
&=&[26,9,12,,4,14,6,14,6,-5,-5,2,2,2,2,\\
&~&\ 4,4,18,10,4,5,7,8,7,8].
\end{array}
\end{align*}
\end{exa}

\section{Inner Product and Orthogonal Decomposition}
Consider any two vectors $X$ and $Y$ in the vector space $\R^{nk}$. The standard inner product on Euclidean space is
$$\langle X,Y\rangle=X^\textrm{T}Y.$$
For positive definite matrix $Q\in {\cal M}_{nk\times nk}$, the weighted inner product is defined as follows
$$\langle X,Y\rangle_Q:=X^\textrm{T}QY,$$
where $Q\in {\cal M}_{nk\times nk}$ is called the weight matrix of the inner product $\langle \cdot,\cdot\rangle_Q$.

Different decompositions employ different inner products. The inner product in \cite{can11} is the weighted inner product with the weight matrix
\begin{align}\label{eq01}
Q=\textrm{diag}\left( \underbrace{k_1,\cdots,k_1}_k,\underbrace{k_2,\cdots,k_2}_k,\cdots,\underbrace{k_n,\cdots,k_n}_k\right).
\end{align}
And the inner product used in \cite{che16} and \cite{shh16} is the standard inner product on Euclidean space. We will investigate under what conditions the two inner products can induce the same decomposition.

Consider a decomposition of finite games, which has the following form
\begin{align}\label{eq1}
{\cal G}_{[n;k_1,\cdots,k_n]}={\cal M}_1~\oplus~{\cal M}_2~\oplus\cdots\oplus~{\cal M}_p,
\end{align}
where ${\cal M}_i$ is the subspace of finite games, $i=1,\ldots,p$.
\begin{dfn} (Compatible Condition)
Consider the decomposition of finite games (\ref{eq1}).  The standard inner product and the weighted inner product with $Q$ as its weight matrix are called compatible, if  for any game  $G\in {\cal M}_{i_j}$ with $V_G$ as its payoff vector, the games determined by $V_GQ$ and $V_GQ^{-1}$ also belong to ${\cal M}_{i_j}, j=1,\cdots,p-1$.
\end{dfn}

\begin{thm}
A common decomposition (\ref{eq1}) can be induced by the standard inner product and a weighted inner product simultaneously, if and only if the compatible condition is satisfied.
\end{thm}
\noindent{\it Proof:}  Suppose the decomposition (\ref{eq1}) is induced by the standard inner product. Then for any $G_1\in {\cal M}_{i_1}$ and $G_2\in {\cal M}_{i_2}$ we have,
\begin{align}\label{eq3}
\begin{array}{ccl}
&~&\langle V_G^{i_1},V_G^{i_2}\rangle\\
&=&V_G^{i_1}(V_G^{i_2})^T\\
&=&0,~\forall i_1,i_2.
\end{array}
\end{align}
Consider a weighted inner product with $Q$ as its weight matrix, we have
\begin{align}\label{eq4}
\begin{array}{ccl}
&~&\langle V_G^{i_1},V_G^{i_2}\rangle_Q\\
&=&V_G^{i_1}Q(V_G^{i_2})^T\\
&=&\hat{V}_G^{i_1}(V_G^{i_2})^T\\
&=&\langle\hat{V}_G^{i_1},(V_G^{i_2})^T\rangle\\
&=&0,~\forall i_1,i_2,
\end{array}
\end{align}
where $\hat{V}_G^{i_1}=V_G^{i_1}Q.$
The first equality follows from the definition of weighted inner product, the second equality follows from the compatible condition, and the third equality follows from condition (\ref{eq3}).
Condition (\ref{eq4}) implies that the decomposition (\ref{eq1}) is induced by the weighted inner product.

If the decomposition (\ref{eq1}) is induced by the weighted inner product. Then for any $G_1\in {\cal M}_{i_1}$ and $G_2\in {\cal M}_{i_2}$ we have,
\begin{align}\label{eq5}
\begin{array}{ccl}
&~&\langle V_G^{i_1},V_G^{i_2}\rangle_Q\\
&=&V_G^{i_1}Q(V_G^{i_2})^T\\
&=&0,~\forall i_1,i_2.
\end{array}
\end{align}
Consider the standard inner product, we have
\begin{align}\label{eq6}
\begin{array}{ccl}
&~&\langle V_G^{i_1},V_G^{i_2}\rangle\\
&=&V_G^{i_1}(V_G^{i_2})^T\\
&=&(V_G^{i_1}Q^{-1})Q(V_G^{i_2})^T\\
&=&\langle V_G^{i_1}Q^{-1},V_G^{i_2}\rangle_Q\\
&=&0,~\forall i_1,i_2.
\end{array}
\end{align}
The first equality follows from the definition of standard inner product, the second and the third equality follows from the compatible condition, and the fourth equality follows from condition (\ref{eq5}).
Condition (\ref{eq6}) implies that the decomposition (\ref{eq1}) is induced by the standard inner product.

\hfill $\Box$
\begin{exa}
\begin{enumerate}
\item Consider the potential based decomposition of finite games
\begin{align*}
{\cal G}_{[n;k_1,\cdots,k_n]}={\cal P}\quad\oplus\quad{\cal N}\quad\oplus\quad{\cal H}.
\end{align*}
According to \cite{che16} and \cite{can11}, the potential based decomposition was obtained using the standard inner product and weighted inner product (\ref{eq01}), respectively. The reason can be explained as follows. For any potential game $G=\{N,S,\{c_i\}_{i\in N}\}$,  the game $\hat{G}$ is a weighted potential game, where $\hat{G}$ is determined by the payoff vector $V_GQ\ (\text{or~} V_GQ^{-1})$. Similarly, it is easy to verify that  any non-strategic game $G$, the game $V_GQ\ (\text{or~} V_GQ^{-1})$ is also non-strategic. Therefore the  compatible condition is satisfied.
\item Consider the  zero-sum based decomposition
\begin{align*}
{\cal G}_{[n;k_1,\cdots,k_n]}={\cal Z}~\oplus~{\cal C}.
\end{align*}
Let $V_G\in {\cal Z}$ and $\tilde{V}_G\in {\cal C}$ be arbitrary. Then
\begin{align*}
\begin{array}{ccl}
&~&\langle V_G,\tilde{V}_G\rangle\\
&=&\sum_{s\in S}\sum_{i=1}^nV_i(s)\tilde{V}_i(s)\\
&=&\sum_{s\in S}\tilde{V}_1(s)\sum_{i=1}^nV_i(s)\\
&=&0,
\end{array}
\end{align*}
which implies that the zero-sum based decomposition can be induced by the  standard inner product. But for weighted inner product (\ref{eq01})
\begin{align}\label{eq7}
\begin{array}{ccl}
\langle V_G,\tilde{V}_G\rangle_Q=\sum_{s\in S}\sum_{i=1}^nk_iV_i(s)\tilde{V}_i(s).
\end{array}
\end{align}
According to (\ref{eq7}), $\langle V_G,\tilde{V}_G\rangle_Q=0$, if and only if all players have the same strategies, i.e. $k_i=constant,\ \forall i.$ Therefore the zero-sum based decomposition does not hold for general case. The reason is that the compatible condition is not satisfied for zero-sum based decomposition.
\item Using the compatible condition, one can verify that the symmetry based decomposition and normalization based decomposition are the same under the standard inner product and  the weighted inner product (\ref{eq01}).
\end{enumerate}
\end{exa}
\section{Concluding Remarks}

This paper considers the effect of different inner products on the  orthogonal decomposition of finite games.  We find that only when the compatible condition is satisfied,  a common decomposition  can be induced by the standard inner product and the weighted inner product simultaneously. To explain the result, we studied the existing decompositions, including potential based decomposition, zero-sum based decomposition, and symmetry based decomposition.

\end{document}